\documentclass[reqno,11pt]{amsart}

\usepackage{amsfonts,amsmath,amssymb,fullpage,graphicx}
\usepackage{verbatim}

\newcommand{\de}{\mathrm{d}}
\newcommand{\e}{\mathrm{e}}
\newcommand{\nn}{\nonumber}
\newcommand{\ug}{\!\!\!\!&=&\!\!\!\!}

\begin{document}

\makeatletter
\title{The algebra of factorial polynomials}
\author{D. Babusci}
\address{INFN - Laboratori Nazionali di Frascati, via  E. Fermi, 40, IT 00044 Frascati (Roma), Italy}
\email{danilo.babusci@lnf.infn.it}
\author{G. Dattoli}
\address{ENEA - Centro Ricerche Frascati, via  E. Fermi, 45, IT 00044 Frascati (Roma), Italy; 
Universit\'e Paris XIII, LIPN, Institut Galil\'ee, CNRS UMR 7030, 99 Av. J.-B. Clement, F 93430 
Villetaneuse, France}
\email{giuseppe.dattoli@enea.it}
\author{M. Carpanese}
\address{ENEA - Centro Ricerche Frascati, via E. Fermi, 45, IT 00044 Frascati (Roma), Italy}
\email{mariano.carpanese@enea.it}

\date{}

\begin{abstract}
We discuss the formal aspects of the factorial polynomials and of the associated series. We develop the theory using 
the formalism of quasi-monomials and prove the usefulness of the method for the solutions of nontrivial difference equations.
\end{abstract}

\maketitle

In this note we consider the so called factorial polynomials \cite{Smirnov}, \cite{Blasiak} by framing the relevant theory within the context 
of quasi-monomials \cite{Dattoli}. The method we propose allows a significant simplification of the study of their properties and of their 
applications in the theory of multi-loop Feynman integrals \cite{Laporta}.

The factorial polynomials (f.p.)
\begin{equation}
\varphi_n (x) = \frac{\Gamma (x + 1)}{\Gamma (x + 1 - n)}
\end{equation} 
are quasi-monomials in the sense that they behave as ordinary monomials under the action of the operators  
\begin{equation}
\label{eq:MP}
\hat{M} = x\,\e^{- \partial_x} \qquad \qquad \hat{P} = \e^{\partial_x} - 1\,,
\end{equation}
namely\footnote{Though not in the explicit realization of Eq. \eqref{eq:MP}, the operators $\hat{M}$ and $\hat{P}$ have been introduced in ref. 
\cite{Dattoli}.}
\begin{equation}
\label{eq:MPact}
\hat{M}\,\varphi_n (x) = \varphi_{n + 1} (x) \qquad\qquad \hat{P}\,\varphi_n (x) = n\,\varphi_{n - 1} (x)\,.
\end{equation}
These operators are called multiplicative and derivative operators, respectively, and satisfy the commutation relation
\begin{equation}
\left[\hat{P}, \hat{M}\right] = 1\,.
\end{equation}

The polynomials $\varphi_n (x)$ satisfy the ``differential" equation
\begin{equation}
\hat{M}\,\hat{P}\,\varphi_n (x) = n\,\varphi_n (x)\,,
\end{equation}
that, according to eq. \eqref{eq:MP}, is equivalent to the difference equation 
\begin{equation}
x\,\left[\varphi_n (x) - \varphi_n (x - 1)\right] = n\,\varphi_n (x)\,. 
\end{equation}

By defining the operator
\begin{equation}
\phi (\hat{M}) = \sum_{n = 0}^\infty a_n\,\hat{M}^n\,,
\end{equation}
the identity 
\begin{equation}
\label{eq:Mn}
\hat{M}^n \,\varphi_0 (x) = \hat{M}^n \,1 =\varphi_n (x)
\end{equation}
allows us to define a f.p.-based function as follows
\begin{equation}
f (x) = \phi (\hat{M})\,1 = \sum_{n = 0}^\infty a_n\,\varphi_n (x)\,.
\end{equation}
For example, to the operator $\exp(\lambda\,\hat{M})$ it is associated the f.p.-based function
\begin{equation}
\label{eq:ex}
e (x, \lambda) = \e^{\lambda\,\hat{M}}\,1 =\sum_{n = 0}^\infty \frac{\lambda^n}{n!}\,\varphi_n (x) = (1 + \lambda)^x
\end{equation}
that plays the same role of exponential function for ordinary monomials since it is an eigenfunction of the operator $\hat{P}$ 
with eigenvalue $\lambda$
\begin{equation}
\hat{P}\,e (x, \lambda) = \lambda\,e (x, \lambda)\,.
\end{equation}
In terms of it, the f.p.-based cosine and sine functions can be defined as follows
\begin{eqnarray}
c (x) \ug \frac12\,\left[e (x, i) + e (x, - i)\right]  = \sum_{k = 0}^\infty \frac{(- 1)^k}{(2\,k)!}\,\hat{M}^{2 k}\,1 = 
\sum_{k = 0}^\infty \frac{(- 1)^k}{(2\,k)!}\,\varphi_{2 k} (x) \nn \\
s (x) \ug \frac1{2\,i}\,\left[e (x, i) - e (x, - i)\right]  = \sum_{k = 0}^\infty \frac{(- 1)^k}{(2\,k + 1)!}\,\hat{M}^{2 k + 1}\,1 = 
\sum_{k = 0}^\infty \frac{(- 1)^k}{(2\,k + 1)!}\,\varphi_{2 k + 1} (x)\,. 
\end{eqnarray}

The f.p.-based function defined in eq. \eqref{eq:ex} is involved in the solution of the finite difference heat equation
\begin{equation}
\label{eq:heat}
\hat{P}_\tau \,W (x, \tau) = \partial_x^2\,W (x,\tau) \qquad\qquad W (x,0) = w (x)
\end{equation}
where
\begin{equation}
\label{eq:difh}
\hat{P}_\tau \,W (x, \tau) = W (x, \tau + 1) - W (x, \tau)\,.
\end{equation}
The solution formally writes\footnote{The suffix $\tau$ indicates that we consider only the space of this variable.}
\begin{equation}
W (x, \tau) = \exp(\hat{M}_\tau\,\partial_x^2)\,w (x)\,1_\tau\,,
\end{equation}
that, taking into account eq. \eqref{eq:ex}, allows us to get the solution of eq. \eqref{eq:heat} in the form
\begin{equation}
W (x, \tau) = e (\tau, \partial_x^2)\, w (x)\,.
\end{equation}
The expansion in terms of factorial series is 
\begin{equation}
W (x, \tau) = \sum_{n = 0}^m \frac{\varphi_n (\tau_m)}{n!}\,w^{(2\,n)} (x)\,,
\end{equation}
where the upper index $(2 n)$ denotes the order of the derivative and $\tau_m = m\,\tau_0$ is the discrete time, corresponding to 
the definition given in eq. \eqref{eq:difh} of the time derivative ($\tau_0 =1$ because we have chosen a unity time interval). The sum 
run up to $m$ as a consequence of the fact that $\varphi_{m + 1} (m) = 0$. The method could be usefully exploited to deal with 
discrete time Schr\"odinger equations \cite{Caldirola}.
\vspace{0.5cm}

In the same way, starting with the operators  
\begin{equation}
\label{eq:Ope}
\beta_n (\hat{M}) = \sum_{k = 0}^\infty \frac{(- 1)^k}{k!\,(n + k)!\,2^{n + 2 k}}\,\hat{M}^{n + 2 k}\,, 
\end{equation}
the f.p.-based Bessel functions can be defined as $B_n (x) = \beta (\hat{M})\,1$. In analogy with the ordinary 
cylindrical Bessel functions, the operators \eqref{eq:Ope} satisfy the recurrences
\begin{eqnarray}
\label{eq:Bes}
2\,\hat{P}\,\beta_n (\hat{M}) \ug \left[\beta_{n - 1} (\hat{M})  - \beta_{n + 1} (\hat{M}) \right]  \nn \\
2\,n\,\beta_n (\hat{M}) \ug \hat{M}\,\left[\beta_{n - 1} (\hat{M})  + \beta_{n + 1} (\hat{M}) \right]\,.
\end{eqnarray}
that, in terms of finite difference, gives
\begin{eqnarray}
B_n (x + 1) \ug B_n (x) + \frac12\,\left[B_{n - 1} (x)  - B_{n + 1} (x) \right]\nn \\
2\,n\,B_n (x) \ug x\,\left[B_{n - 1} (x - 1)  + B_{n + 1} (x - 1) \right]\,.
\end{eqnarray}
By denoting with $\hat{N}$ a number operator such that $\hat{N}\,B_n (x) = n \,B_n (x)$, the combination of 
eqs. \eqref{eq:Bes} allows the introduction of the following index shifting operators
\begin{equation}
\hat{E}_+ = - \hat{P} + \hat{M}^{- 1}\,\hat{N}\,, \qquad \qquad \hat{E}_- = \hat{P} + \hat{M}^{- 1}\,\hat{N}\,,
\end{equation}
whose action on the operators in eq. \eqref{eq:Ope} is given by
\begin{equation}
\hat{E}_\pm\,\beta_n (\hat{M}) = \beta_{n \pm 1} (\hat{M})
\end{equation}
from which it is easy to show that 
\begin{equation}
\left[(\hat{M}\,\hat{P})^2 + (\hat{M}^2 - n^2)\right]\,\beta_n (\hat{M}) = 0\,.
\end{equation}
The corresponding difference equation writes
\begin{equation}
2\,x\,(x - 1)\,B_n (x - 2) - x\,(2\,x - 1)\,B_n (x - 1) + (x^2 - n^2)\,B_n (x) = 0\,.
\end{equation}
\vspace{0.5cm}

So far a differential equation has been translated into a corresponding difference equation using the realization \eqref{eq:MP} 
of the multiplicative and derivative operators. The correspondence between the solutions (i.e. their isospectrality) is ensured by the 
fact that the operators $\hat{P}$, $\hat{M}$, and $\hat{1}$ realize a Weyl algebra. Let us now consider the same problem from a 
reversed perspective, i.e., given a difference equation we obtain from it a corresponding differential equation in terms of derivative and 
multiplicative operators.

Let us consider the following difference equation
\begin{equation}
\label{eq:diff}
(4\,x + 2)\,F (x + 1) + 4\,x\,F(x - 1) - (8\,x - 1)\,F (x) = 0\,.
\end{equation}
By performing the following substitutions
\begin{equation}
\e^{\partial_x} \to 1 + \hat{P}\,, \qquad\qquad x \to \hat{M}\,(1 + \hat{P})\,,
\end{equation}
we find 
\begin{equation}
\left(4\,\hat{M}\,\hat{P}^2 + 2\,\,\hat{P} + 1\right)\,\Phi (\hat{M}) = 0
\end{equation}
where $\Phi (\hat{M}))\,1 = F (x)$. The solution of the previous differential equation can be found by using the Frobenius method. 
Hence, we look for solutions of the type
\begin{equation}
\Phi (\hat{M}) = \sum_{k = 0}^\infty a_k\,M^{k + c}
\end{equation}
and the associated indicial equation is
\begin{equation}
2\,a_0\,c\,(2\,c - 1) = 0\,.
\end{equation}
By using the solution $c = 0$ of this equation, we get
\begin{equation}
F (x) = \Phi (\hat{M})\,1 = \sum_{k = 0}^\infty \frac{( - 1)^k}{(2\,k)!}\,\hat{M}^k = \sum_{k = 0}^\infty \frac{( - 1)^k}{(2\,k)!}\,\varphi_k (x)\,.
\end{equation}
The other solution, $c = 1/2$, implies the definition of the operator $\hat{M}^{1/2}$, whose meaning will be discussed in a forthcoming 
paper, where the problems associated with the fractional powers of operators of the type \eqref{eq:MP} will be considered.

The method we have envisaged can also be extended to the non-homogeneous case. For example, the equation
\begin{equation}
a\,y (x + 1) + b\,y (x) = g (x)
\end{equation}
can be cast in the form
\begin{equation}
\hat{Q}\, \psi (\hat{M}) = g (\hat{M}\,(1 + \hat{P}))
\end{equation}
where $\hat{Q} = a\,\hat{P} + (a + b)$ and $y (x) = \psi (\hat{M})\,1$. The particular solution of this equation formally writes
\begin{equation}
\psi_p (\hat{M}) = \hat{Q}^{- 1}\,g (\hat{M}\,(1 + \hat{P}))
\end{equation}
and the explicit form can be written as
\begin{equation}
\psi_p (\hat{M}) = \int_0^\infty \de s\,\e^{- s\,\hat{Q}}\,g (\hat{M}\,(1 + \hat{P}))
\end{equation}
from which, taking into account that
\begin{equation}
\e^{\alpha\,\hat{P}}\,g (\hat{M}\,(1 + \hat{P})) = g ((\hat{M} + \alpha)\,\,(1 + \hat{P}))\,,
\end{equation}
we obtain
\begin{equation}
\psi_p (\hat{M}) = \int_0^\infty \de s\,\e^{- s\,(a + b)}\,g ((\hat{M} - a\,s)\,(1 + \hat{P}))\,.
\end{equation}
As an example, in the case $g (x) = x^2$, by using eqs. \eqref{eq:MPact} and \eqref{eq:Mn}, one obtains:
\begin{equation}
y_p (x) = \frac{z}{a}\,\left\{(x - z)^2 - z\,(1 - z)\right\} \qquad\qquad \left(z = \frac{a}{a + b}\right)\,.
\end{equation}

This note has just been aimed at fixing the general rule of the use of factorial series in the theory of difference equations. More 
specific problems regarding applications like those treated in Ref. \cite{Laporta} will be discussed elsewhere.

\section*{Acknowledgments}
One of the authors (G. D.) has benefited from interesting and enlightening discussions with Prof. G. H. E. Duchamp at University of 
Paris XIII, where this work was started. He also recognizes the warm hospitality and the financial support of 
Universit Paris XIII, LIPN, Institut GalilŽe, CNRS UMR 7030.

\end{document}